# On approximate pseudo-maximum likelihood estimation for LARCH-processes

JAN BERAN and MARTIN SCHÜTZNER

*Department of Mathematics and Statistics, University of Konstanz, Universitätsstrasse 10, 78457 Konstanz, Germany. E-mail: jan.beran@uni-konstanz.de; martin.schuetzner@uni-konstanz.de*

Linear ARCH (LARCH) processes were introduced by Robinson [*J. Econometrics* **47** (1991) 67–84] to model long-range dependence in volatility and leverage. Basic theoretical properties of LARCH processes have been investigated in the recent literature. However, there is a lack of estimation methods and corresponding asymptotic theory. In this paper, we consider estimation of the dependence parameters for LARCH processes with non-summable hyperbolically decaying coefficients. Asymptotic limit theorems are derived. A central limit theorem with $\sqrt{n}$-rate of convergence holds for an approximate conditional pseudo-maximum likelihood estimator. To obtain a computable version that includes observed values only, a further approximation is required. The computable estimator is again asymptotically normal, however with a rate of convergence that is slower than $\sqrt{n}$.

## 1. Introduction

Since the introduction of ARCH and GARCH processes in the seminal papers of Engle (1982) and Bollerslev (1986), an abundance of models with conditional heteroskedasticity have been proposed. More recently, modifications of these models have been introduced to include the possibility of slowly decaying correlations (long memory) in volatility. This was motivated by the observation that empirical autocorrelations in squared log-returns often persist over long stretches of time. Long memory means that the sum of autocorrelations over all lags is infinite. As it turns out, not all models proposed in this context have long memory in volatility, although their correlations may decay hyperbolically. For instance, no second order stationary ARCH($\infty$) process $X_t$ with non-summable autocorrelations of $X_t^2$ exists (Giraitis *et al.* (2000a, 2000b)). Models with genuine long memory in volatility include linear ARCH (LARCH) models introduced by Robinson (1991) and stochastic volatility (SV) models such as the FIEGARCH process (Harvey (1998), Robinson (2001), Surgailis and Viano (2002)). With respect to estimation, SV models are somewhat complicated since they are based on unobservable latent processes. In contrast, no latent process is included in the definition of LARCH processes. This allows for direct estimation of unknown parameters, including maximum likelihood estimation and related methods. For LARCH processes, the difficulty in studying asymptotics of parameter estimates is, however, the rather complex structure of the stationary solution (Giraitis *et al.* (2000a, 2000b)). The problem of location estimation is considered in Beran (2006). Related limit theorems can be found in Berkes and Horvath (2003) and Giraitis *et al.* (2000a, 2000b). Here, we will consider estimation of dependence parameters for LARCH processes with hyperbolically decaying non-summable weights.

A LARCH process $(X_t, \sigma_t)_{t \in \mathbb{Z}}$ is defined by

$$X_t = \varepsilon_t \sigma_t, \tag{1}$$

$$\sigma_t = a + \sum_{j=1}^{\infty} b_j X_{t-j}, \tag{2}$$

where the following assumptions hold:

(A1) $\varepsilon_t$ are i.i.d. random variables defined on a probability space $(\Omega, \mathcal{A}, P)$, with continuous distribution, $E(\varepsilon_t) = 0$, and $E(\varepsilon_t^2) = 1$;
(A2) $a \neq 0$ and $b = \sum_{j=1}^{\infty} b_j^2 < 1$.

The stationary solution of the LARCH equations is given by

$$\sigma_t = a + a \sum_{k=1}^{\infty} \sum_{j_1,\ldots,j_k=1}^{\infty} b_{j_1} \cdots b_{j_k} \varepsilon_{t-j_1} \cdots \varepsilon_{t-j_1-\cdots-j_k}$$

(Giraitis *et al.* (2000a, 2000b)). Obviously, the process $(X_t)_{t \in \mathbb{Z}}$ is uncorrelated. Giraitis *et al.* (2003) showed that if $b_j \sim_{j \to \infty} c j^{d-1}$ for some $d \in (0, \frac{1}{2})$ and $E(X_t^4) < \infty$, then there is long memory in volatility characterized by

$$\gamma_\sigma(k) = \text{cov}(\sigma_0, \sigma_k) \underset{|k| \to \infty}{\sim} c_1 |k|^{2d-1}$$

and

$$\gamma_{X^2}(k) = \text{cov}(X_0^2, X_k^2) \underset{|k| \to \infty}{\sim} c_2 |k|^{2d-1},$$

and the same is true for the leverage covariance $\gamma_L(k) = \text{cov}(\sigma_k^2, X_0)$.

The main purpose of our work is to provide statistical theory for the estimation of a parametric version of (1) and (2). Thus, we assume $a$ and $(b_j)_{j \geq 1}$ to depend on a finite-dimensional parameter vector $\theta$. We will focus on conditional maximum likelihood estimation, a method often used for models with conditional heteroskedasticity. Under the assumption of Gaussian $\varepsilon_t$, the following approximate maximum likelihood estimator of $\theta$ can be defined:

$$\theta_n^* := \arg\min_{\theta \in \Theta} L_n^*(\theta),$$

where

$$L_n^*(\theta) = \sum_{t=1}^n \frac{X_t^2}{\sigma_t^2(\theta)} + \ln \sigma_t^2(\theta)$$

and

$$\sigma_t(\theta) = a(\theta) + \sum_{j=1}^\infty b_j(\theta) X_{t-j}.$$

Given a finite sample, $\sigma_t(\theta)$ has to be replaced by a proxy $\bar{\sigma}_t(\theta)$, depending on the finite past only. Since, in general, $\varepsilon_t$ is not assumed to be normal, $\theta_n^*$ is called a pseudo-maximum likelihood estimator (PMLE). In the case where $(X_t, \sigma_t)$ is the original ARCH(1) or GARCH(1,1) process, the asymptotic properties of $\theta_n^*$ have been investigated in Lee and Hanson (1994) and Lumsdaine (1996), and were generalized to GARCH($p,q$) and ARCH($\infty$) processes by Berkes *et al.* (2003) and Robinson and Zaffaroni (2006), respectively. For long-memory LARCH processes, derivation of asymptotic results is more complicated because the coefficients $b_j$ are not summable. Moreover, $\sigma_t^2$ may become arbitrarily small and hence $\sigma_t^{-2}$ and its derivatives arbitrarily large. The first problem leads to difficulties with respect to differentiability of $\sigma_t(\theta)$ as a function of $\theta$. Additional assumptions on the parametric model are therefore needed (see Section 2). The second problem can be avoided by modifying the original maximum likelihood equations (see Section 3). Also, note that parametric estimation for finite order LARCH processes, that is, where the sum in (2) is finite and thus the autocorrelations of the squares are absolutely summable, is considered in Francq and Zakoian (2008) and Truquet (2008).

The outline of the paper is as follows. Section 2 deals with ergodicity and differentiability as necessary prerequisites. Estimation of $\theta$ is considered in Section 3. Asymptotic results are derived for two versions of a modified MLE: (a) estimate with $\sigma_t(\theta)$ ($t=1,\ldots,n$) and (b) estimate including only values of $\sigma_t(\theta)$ that can be approximated with sufficient accuracy. Lemmas needed in the proofs of the main results can be found in the Appendix. A small simulation study in Section 4 illustrates the theoretical results. Some general comments in Section 5 conclude the paper.

## 2. Ergodicity and differentiability

### 2.1. Ergodicity

To ensure consistency, ergodicity of $\sigma_t$ is needed. The following proposition is an extension of Theorem 2.1 in Giraitis *et al.* (2003).

**Proposition 1.** *Under* (A1) *and* (A2), *there exists a unique strictly and second order stationary solution of (1) and (2). This solution is ergodic.*

**Proof.** $\sigma_t$ is given by the Volterra decomposition (see Giraitis *et al.* (2000a, 2000b))

$$\sigma_t = a + a \sum_{k=1}^{\infty} \sum_{j_1,\ldots,j_k=1}^{\infty} b_{j_1} \cdots b_{j_k} \varepsilon_{t-j_1} \cdots \varepsilon_{t-j_1-\cdots-j_k}.$$

Since $\{\varepsilon_{i_1} \cdots \varepsilon_{i_r}\}_{1 \leq i_1 < \cdots < i_r, r \geq 1}$ is an orthonormal system, convergence in the $L^2(\Omega)$-norm follows from (A2) since

$$\sum_{k=1}^{\infty} \sum_{j_1,\ldots,j_k=1}^{\infty} b_{j_1}^2 \cdots b_{j_k}^2 = \sum_{k=1}^{\infty} b^k < \infty.$$

For the uniqueness of $\sigma_t$, we refer to Giraitis *et al.* (2003). For the proof of ergodicity, it is sufficient to find a measurable function $f : \mathbb{R}^{\infty} \to \mathbb{R}$ with $\sigma_t = f(\varepsilon_{t-1}, \varepsilon_{t-2}, \ldots)$, where equality holds *almost surely* (see, for example, Theorem 3.5.8 in Stout (1974)). First, note that convergence of the infinite sum defining the solution is independent of the order of summation since the series of squared coefficients is absolutely summable. Hence, we make use of the following alternative representation of $\sigma_t$. Define

$$f_k(x_1, x_2, \ldots) = \sum_{\substack{j_i \geq 1, l \leq k \\ j_1 + \cdots + j_l = k}} b_{j_1} \cdots b_{j_l} x_{j_1} \cdots x_{j_1 + \cdots + j_l}$$

and

$$M_t(k) = f_k(\varepsilon_{t-1}, \varepsilon_{t-2}, \ldots).$$

Then

$$\sigma_t = a + a \sum_{k=1}^{\infty} M_t(k)$$

and for every fixed $t \in \mathbb{Z}$, $M_t(k), k = 1, 2, \ldots$, is a martingale difference w.r.t. $\mathcal{F}_k^t = \sigma\{M_t(l), l \leq k\}$. An application of the martingale convergence theorem yields that

$$S_t(m) = \sum_{k=1}^{m} M_t(k) \to \sum_{k=1}^{\infty} M_t(k)$$

as $m \to \infty$ almost surely. Hence, the desired representation is given by

$$f = \sum_{k=1}^{\infty} f_k.$$

For the measurability of $f$, see Corollary 2.1.3 in Straumann (2004). $\square$

## 2.2. Differentiability

For simplicity of notation, we will concentrate on coefficients $(b_j)_{j\geq 1}$ of the following type:

(B1)
$$b_j(c,d) = cj^{d-1},$$

where $d \in [0, d_u]$, $d_u < \frac{1}{2}$, $c \in [0, c_u(d)]$ and

$$c_u(d) = C\left(\sum_{j=1}^{\infty} j^{2d-2}\right)^{-1/2},$$

with $0 < C < 1$.

(B2) $a \in [a_d, a_u]$ with $0 < a_d < a_u < \infty$.

Assumption (B1) ensures the summability constraint in (A2). Extending the results to more general weights, such as, for instance, those obtained from the FARIMA$(p,d,q)$ operator (see Grager and Joyeux (1980), Hosking (1981)) is straightforward. For instance, we may consider FARIMA$(0,d,0)$ weights $b_j$ defined by

$$\sum_{j=1}^{\infty} b_j B^j = c(d)[(1-B)^{-d} - 1],$$

where $0 < d < \frac{1}{2}$ and $c(d)$ is a constant such that $\sum b_j^2 < 1$. Note, in particular, that here $(1-B)^{-d}$ instead of $(1-B)^d$ induces long memory for $d > 0$.

In the following, we will use the notation $\Theta \subset [0, \frac{1}{2}) \times (\mathbb{R}^+)^2$ for the set of all $\theta = (d, c, a)^T$ such that (B1) and (B2) hold. Moreover, for a real matrix $A$, we define the matrix norm

$$\|A\| = \operatorname{tr}(A^T A)^{1/2}.$$

Convergence of matrices will be understood with respect to this norm. The LARCH process $(X_t, \sigma_t)_{t \in \mathbb{Z}}$ will be assumed to belong to the parametric family with $\theta_0$ in the interior of $\Theta$.

From the given dynamical structure in equation (2), we can reconstruct the unobservable conditional variance $\sigma_t^2$ from the infinite past $(X_s)_{s \leq t}$, as follows. Define, for any $\theta \in \Theta$ and $t \in \mathbb{Z}$,

$$\sigma_t(\theta) = a + \sum_{j=1}^{\infty} b_j(c,d) X_{t-j}.$$

For the process with true parameter $\theta_0$, we have, in particular,

$$\sigma_t^2(\theta_0) = \operatorname{var}(X_t \mid X_s, s \leq t-1).$$

Given a finite sample $(X_t)_{t=1,\ldots,n}$, $\sigma_t(\theta)$ has to be approximated, for instance, by

$$\bar{\sigma}_t(\theta) = a + \sum_{j=1}^{t-1} b_j(c,d) X_{t-j}, \qquad t \geq 1.$$

The extent to which this may be a good approximation of $\sigma_t(\theta)$ will be discussed in Section 3.

We now consider the properties of $\sigma_t(\theta)$ for fixed $t \in \mathbb{Z}$ as a stochastic process with index $\theta \in \Theta$. The reason is that almost sure continuity and differentiability of $\sigma_t(\theta)$ as a function of $\theta$ will be required in the next section. Moreover, we need to ensure measurability of infima involving $\sigma_t(\theta)$ on the uncountable set $\Theta$. In the case of absolutely summable coefficients $(b_j)_{j\geq 1}$, this is not a problem since the infinite sum defining the stationary solution is uniformly absolutely summable, on a set of probability one, and $\sigma_t(\theta)$ inherits the properties of $b_j(c,d)$. In contrast, for non-summable $b_j$, this is not automatically the case. We therefore impose the following assumption.

(S) For every $t \in \mathbb{Z}$, $(\sigma_t(\theta))_{\theta \in \Theta}$ is a separable stochastic process on $\Theta$, that is, for every open $A \subset \Theta$ and closed interval $B$, the sets

$$\{\omega | \sigma_t(\theta) \in B, \forall \theta \in A\} \quad \text{and} \quad \{\omega | \sigma_t(\theta) \in B, \forall \theta \in A \cap \mathbb{Q}^3\}$$

differ only on a set $N \subset N_0$, where $P(N_0) = 0$.

**Remark 1.** The process $(\sigma_t(\theta))_{\theta \in \Theta}$ can always be replaced by a separable version (see Theorem 2.4 in Doob (1953)).

The following result can now be obtained.

**Proposition 2.** *Under assumptions* (A1), (B1), (B2) *and* (S), $\sigma_t(\theta)$ *is almost surely infinitely often differentiable in $\theta$ and the kth partial derivative w.r.t. d is given by*

$$\frac{\partial^k}{\partial d^k} \sigma_t(\theta) = \sum_{j=1}^{\infty} \frac{\partial^k}{\partial d^k} b_j(c,d) X_{t-j}.$$

**Proof.** Let

$$\sigma_t(d) := \sigma_t\{(1,1,d)^T\}.$$

The covariance function of $[\sigma_t(d)]_{0 \leq d \leq d_u}$ is given by

$$v(d,d') = \text{Cov}(\sigma_t(d), \sigma_t(d')) = \sum_{j=1}^{\infty} j^{d+d'-2},$$

which is infinitely often differentiable for all $0 \leq d, d' < \frac{1}{2}$. Since $a$ and $c$ are just additive and multiplicative components, respectively, in $\sigma_t(\theta)$, existence of derivatives follows immediately from Lemma 1 (see the Appendix). Indeed, iteration of the following calculation shows that the partial derivatives w.r.t. $d$ can be calculated as claimed: Taylor series expansion of $b_j(1,d)$ for each $j$ up to order 2 yields

$$E\left|\frac{1}{h}(\sigma_t(d+h)-\sigma_t(d))-\sum_{j=1}^\infty \frac{\partial}{\partial d}b_j(1,d)X_{t-j}\right|^2 = h^2 E\left|\sum_{j=1}^\infty \frac{\partial^2}{\partial d^2}b_j(1,\tilde{d}_j)X_{t-j}\right|^2$$

$$= h^2 E(\sigma_t^2)\sum_{j=1}^\infty \left(\frac{\partial^2}{\partial d^2}b_j(1,\tilde{d}_j)\right)^2 \to 0$$

as $h \to 0$, where $d \leq \tilde{d}_j \leq d_u$. $\square$

Lemma 1 in the Appendix also implies that, under (S), we are able to find bounds on $E(\sup_{\theta \in \Theta}|\sigma_t(\theta)|^m)$ ($m \geq 1$) in terms of $\sup_{\theta \in \Theta} E(|\sigma_t(\theta)|^m)$ and $\sup_{\theta \in \Theta} E(|\frac{\partial}{\partial \theta}\sigma_t(\theta)|^m)$. This is very useful for proving uniform convergence results.

## 3. Estimation

### 3.1. Estimation with exact conditional variances

Define

$$\mu_p = E(\varepsilon_t^p),$$
$$|\mu|_p = E(|\varepsilon_t|^p)$$

and

$$\|b\|_p^p = \sum_{j=1}^\infty |b_j|^p.$$

The following assumptions ensure the existence of unconditional moments of $\sigma_t$ and $X_t$. Assumptions (M$_3$), (M$_4$) and (M$''_p$) are from Giraitis *et al.* (2003), while (M$'_p$) is from Giraitis *et al.* (2000b).

(M$_3$) $|\mu|_3 < \infty$ and $|\mu|_3^{1/3}\|b(\theta_0)\|_3 + 3\zeta\|b(\theta_0)\|_2 < 1$, where $\zeta$ is the positive solution of the equation $3\zeta^2 - 3\zeta - 1 = 0$.

(M$'_p$) For $p \geq 2$, $|\mu|_p < \infty$ and $(2^p - p - 1)^{1/2}|\mu|_p^{1/p}\|b(\theta_0)\|_2 < 1$.

(M$''_p$) For even $p \geq 4$, $|\mu|_p < \infty$ and $\sum_{j=2}^p \binom{p}{j}\|b(\theta_0)\|_j^j|\mu_j| < 1$.

**Remark 2.** For even $p \geq 4$, (M$''_p$) is weaker than (M$'_p$). For Gaussian (and similar) $\varepsilon_t$, (M$_3$) is weaker than (M$'_3$). We will therefore make use of assumption (M$'_p$) only if either $p = 5$ or (M$'_3$) is weaker than (M$_3$). The assumptions we will use are only sufficient; more general (but complicated) conditions can be formulated in terms of the moments of $\sigma_t(\theta)$ and its derivatives.

First, we will assume that $\sigma_t(\theta)$ can be calculated exactly, that is, as if we knew the infinite past $(X_s)_{s \leq n}$. To avoid the problem of unbounded $\sigma_t^{-2}$ (see Section 1), we modify the maximum likelihood estimator as follows. Let

$$L_n(\theta) = \frac{1}{n} \sum_{t=1}^{n} l_t(\theta) = \frac{1}{n} \sum_{t=1}^{n} \frac{X_t^2 + \epsilon}{\sigma_t^2(\theta) + \epsilon} + \ln(\sigma_t^2(\theta) + \epsilon),$$

where $\epsilon > 0$ is a small but positive constant, and define the estimator

$$\theta_n^{(1)} := \arg\min_{\theta \in \Theta} L_n(\theta).$$

Furthermore, denote by $L(\theta) = E[l_t(\theta)]$ the expected value of the individual terms in $L_n$. Consistency is given by the following result.

**Theorem 1.** *Let $\epsilon > 0$ and assume that* (A1), (B1), (B2) *and* (S) *hold. Then, under* $(M_3)$ *or* $(M_3')$, $\theta_n^{(1)}$ *is a strongly consistent estimator of $\theta_0$, that is, as $n \to \infty$,*

$$\theta_n^{(1)} \to \theta_0 \qquad a.s.$$

**Proof.** From Lemmas 3 and 4 (see the Appendix), we get uniform a.s. convergence of $L_n(\theta)$ to the function $L(\theta)$. Moreover, $L(\theta)$ has a unique minimum at $\theta_0$. The proof then follows from standard arguments (see, for example, Huber (1967)). □

The asymptotic distribution of $\theta_n^{(1)}$ is essentially determined by $L_n'(\theta_0)$, where

$$L_n'(\theta) = \frac{1}{n} \sum_{t=1}^{n} \frac{\partial}{\partial \theta} l_t(\theta) = \frac{1}{n} \sum_{t=1}^{n} \left(1 - \frac{X_t^2 + \epsilon}{\sigma_t^2(\theta) + \epsilon}\right) \frac{2\sigma_t(\theta)}{\sigma_t^2(\theta) + \epsilon} \frac{\partial}{\partial \theta} \sigma_t(\theta).$$

Define the matrices

$$G_\epsilon = E\left(\frac{\partial}{\partial \theta} l_t(\theta_0) \left(\frac{\partial}{\partial \theta} l_t(\theta_0)\right)^T\right) = E\left(\frac{\sigma_t^4(E\varepsilon_t^4 - 1)}{(\sigma_t^2 + \epsilon)^2} \frac{4\sigma_t^2}{(\sigma_t^2 + \epsilon)^2} \dot\sigma_t \dot\sigma_t^T\right),$$

$$H_\epsilon = E\left(\frac{\partial^2}{\partial \theta \, \partial \theta'} l_t(\theta_0)\right) = E\left(\frac{4\sigma_t^2}{(\sigma_t^2 + \epsilon)^2} \dot\sigma_t \dot\sigma_t^T\right),$$

where

$$\dot\sigma_t = \frac{\partial}{\partial \theta} \sigma_t(\theta_0).$$

The Hessian matrix $\frac{\partial^2}{\partial \theta \, \partial \theta'} l_t(\theta)$ is given explicitly in the proof of Lemma 3 in the Appendix. The asymptotic distribution of $\theta_n^{(1)}$ can now be derived as follows.

**Theorem 2.** *Let $\epsilon > 0$ and $\theta_0$ be in the interior of $\Theta$. Then, under assumptions* (A1), (B1), (B2), (S), (M$'_5$),

$$n^{1/2}(\theta_n^{(1)} - \theta_0) \xrightarrow{d} N(0, H_\epsilon^{-1} G_\epsilon H_\epsilon^{-1})$$

*as $n \to \infty$, where $N(0, \Sigma)$ denotes the three-dimensional centered normal distribution with covariance matrix $\Sigma$.*

**Proof.** By Taylor series expansion,

$$0 = L'_n(\theta_n^{(1)}) = L'_n(\theta_0) + \tilde{L}''_n \cdot (\theta_n^{(1)} - \theta_0)$$

with

$$\tilde{L}''_n = \frac{1}{n} \sum_{t=1}^{n} \frac{\partial^2}{\partial \theta \, \partial \theta'} l_t(\theta),$$

evaluated in each row $j = 1, 2, 3$ at some point $\theta = \tilde{\theta}_n^j$ with $\|\tilde{\theta}_n^j - \theta_0\| \leq \|\theta_n^{(1)} - \theta_0\|$. Since

$$E\left(\frac{\partial}{\partial \theta} l_t(\theta_0) \Big| \mathcal{F}_{t-1}\right) = 0,$$

where $\mathcal{F}_t = \sigma(\varepsilon_s, s \leq t)$, $\frac{\partial}{\partial \theta} l_t(\theta_0)$ is a vector of stationary, ergodic martingale differences with finite variance. Hence, from Theorem 23.1 in Billingsley (1968) and the Cramér–Wold device,

$$n^{1/2} L'_n(\theta_0) \xrightarrow{d} N(0, G_\epsilon)$$

as $n \to \infty$. From Lemma 3 and Proposition 1, we get

$$\tilde{L}''_n \to H_\epsilon$$

almost surely as $n \to \infty$. By Lemma 5, $H_\epsilon$ is invertible. This, together with Slutsky's theorem, concludes the proof. □

*Remark 3.* Letting $\epsilon$ tend to zero, we get $H_\epsilon^{-1} G_\epsilon H_\epsilon^{-1} \to (E\varepsilon_t^4 - 1) H_0^{-1}$, where $H_0 = 4E(\frac{\dot{\sigma}_t \dot{\sigma}_t^T}{\sigma_t^2})$. If $E(\sigma_t^{-2}) = \infty$, this means, for instance, that the asymptotic variance of $\hat{a}$ approaches zero.

*Remark 4.* Formally, we get the same rate of convergence and asymptotic variance as for short memory models, such as GARCH$(p,q)$ and ARCH$(\infty)$ (see Berkes *et al.* (2003) and Robinson and Zaffaroni (2006)).

### 3.2. Estimation given the finite past

Given a finite sample $X_1, \ldots, X_n$, the computable version of the estimator is defined by

$$\theta_n^{(2)} := \arg\min_{\theta \in \Theta} \bar{L}_n(\theta),$$

where

$$\bar{L}_n(\theta) := \frac{1}{n} \sum_{t=1}^{n} \frac{X_t^2 + \epsilon}{\bar{\sigma}_t^2(\theta) + \epsilon} + \ln(\bar{\sigma}_t^2(\theta) + \epsilon).$$

This estimator is consistent in the following sense.

**Theorem 3.** *Let $\epsilon > 0$ and assume that* (A1), (B1), (B2) *and* (S) *hold. Then, under* $(M_3)$ *or* $(M_3')$,

$$\theta_n^{(2)} \to \theta_0$$

*as $n \to \infty$, where convergence holds in $L^1$ and in probability.*

**Proof.** The proof follows as for Theorem 1, with the additional application of Lemma 6. □

Obtaining the asymptotic distribution of $\theta_n^{(2)}$ is more complicated due to the slow convergence of $|\sigma_t(\theta) - \bar{\sigma}_t(\theta)|$ to zero. To be more specific, note that

$$E[(\sigma_t(\theta) - \bar{\sigma}_t(\theta))^2] = \sum_{j=t}^{\infty} b_j^2(c,d) \sim c_1 t^{2d-1}.$$

As in the proof of Theorem 2, Taylor series expansion yields

$$0 = \bar{L}_n'(\theta_n^{(2)}) = \bar{L}_n'(\theta_0) + \widetilde{\bar{L}_n''} \cdot (\theta_n^{(2)} - \theta_0),$$

where $\bar{L}_n'(\theta)$ and $\bar{L}_n''(\theta)$ are the same as $L_n'(\theta)$ and $L_n''(\theta)$ with $\sigma_t(\theta)$ replaced by $\bar{\sigma}_t(\theta)$. Since the law of large numbers still holds (see Lemma 6), the asymptotic distribution of $\theta_n^{(2)}$ follows from the asymptotic distribution of $\bar{L}_n'(\theta_0)$. The latter is the same as for $L_n'(\theta_0)$, provided that

$$d_n := \sqrt{n}(L_n'(\theta_0) - \bar{L}_n'(\theta_0)) \xrightarrow{p} 0$$

as $n \to \infty$. Since $d_n$ is asymptotically equivalent to

$$\frac{1}{\sqrt{n}} \sum_{t=1}^{n} \frac{\dot{\bar{\sigma}}_t(\theta) \bar{\sigma}_t(\theta)(X_t^2 + \epsilon)}{\bar{\sigma}_t^2(\theta) + \epsilon} \left( \frac{1}{\bar{\sigma}_t^2(\theta) + \epsilon} - \frac{1}{\sigma_t^2(\theta) + \epsilon} \right),$$

applying the mean value theorem to $(x^2 + \epsilon)^{-1}$ and taking into account the asymptotic behavior of $E[(\sigma_t(\theta) - \bar{\sigma}_t(\theta))^2]$, a rough upper bound for $E(|d_n|)$ is given by

$$c_1 E\left(\frac{1}{\sqrt{n}} \sum_{t=1}^n \left|\frac{\dot{\bar{\sigma}}_t(\theta)\bar{\sigma}_t(\theta)(X_t^2 + \epsilon)}{\bar{\sigma}_t^2(\theta) + \epsilon}\right| |\sigma_t(\theta) - \bar{\sigma}_t(\theta)|\right) \leq c_2 \frac{1}{\sqrt{n}} n n^{d-1/2}.$$

In the long-memory case with $d > 0$, this bound does not converge to zero. We therefore propose an alternative estimator, at the cost of a slower rate of convergence: for given $0 < \beta < 1$, define $m(n) = \lfloor n^\beta \rfloor - 1$, where $\lfloor \cdot \rfloor$ denotes the floor function,

$$\tilde{L}_n(\theta) := \frac{1}{m(n)} \sum_{t=n-m(n)}^n \frac{X_t^2 + \epsilon}{\bar{\sigma}_t^2(\theta) + \epsilon} + \ln(\bar{\sigma}_t^2(\theta) + \epsilon)$$

and

$$\theta_n^{(\beta)} := \arg\min_{\theta \in \Theta} \tilde{L}_n(\theta).$$

This estimator has the following properties.

**Theorem 4.** *Let $\epsilon > 0$, $\theta_0$ be in the interior of $\Theta$ and assume* (A1), (B1), (B2) *and* (S). *The following then hold:*

(a) *if* (M$_3$) *or* (M$_3'$) *holds and $0 < \beta < 1$, then $\theta_n^{(\beta)}$ converges in $L^1$ and in probability to $\theta_0$;*
(b) *if* (M$_5'$) *holds and $0 < \beta < 1 - 2d$, then as $n \to \infty$,*

$$n^{\beta/2}(\theta_n^{(\beta)} - \theta_0) \xrightarrow{d} N(0, H_\epsilon^{-1} G_\epsilon H_\epsilon^{-1});$$

(c) *if* (M$_3$) *or* (M$_3'$) *holds and $\beta = 1 - 2d$, then*

$$E[|\theta_n^{(\beta)} - \theta_0|] \sim c_2 n^{-(1/2-d)}.$$

**Proof.** The proof is a combination of Theorem 2 and the arguments given above.  □

*Remark 5.* The choice of $\epsilon$ is important for a good performance of the estimator $\theta_n^{(\beta)}$. While the above theorems indicate that $\epsilon$ should be chosen as small as possible, the optimization in the definition of $\theta_n^{(\beta)}$ becomes numerically more demanding if $\epsilon \to 0$ since the function $\tilde{L}_n$ may then exhibit many local minima. As an illustration, in Figure 1, $\tilde{L}_n$ is plotted as a function of the single parameter $d$ for different values of $\epsilon$. How this effect can be handled statistically and how it depends on the parameter $\theta_0$ are the subjects of current research.

*Remark 6.* Calculations analogous to those above imply that for short-memory LARCH processes (that is, LARCH processes with absolutely summable autocorrelations of $X_t^2$), the central limit theorem for $\theta_n^{(2)}$ holds with $\sqrt{n}$-rate of convergence. This also includes the case where $d < 0$.

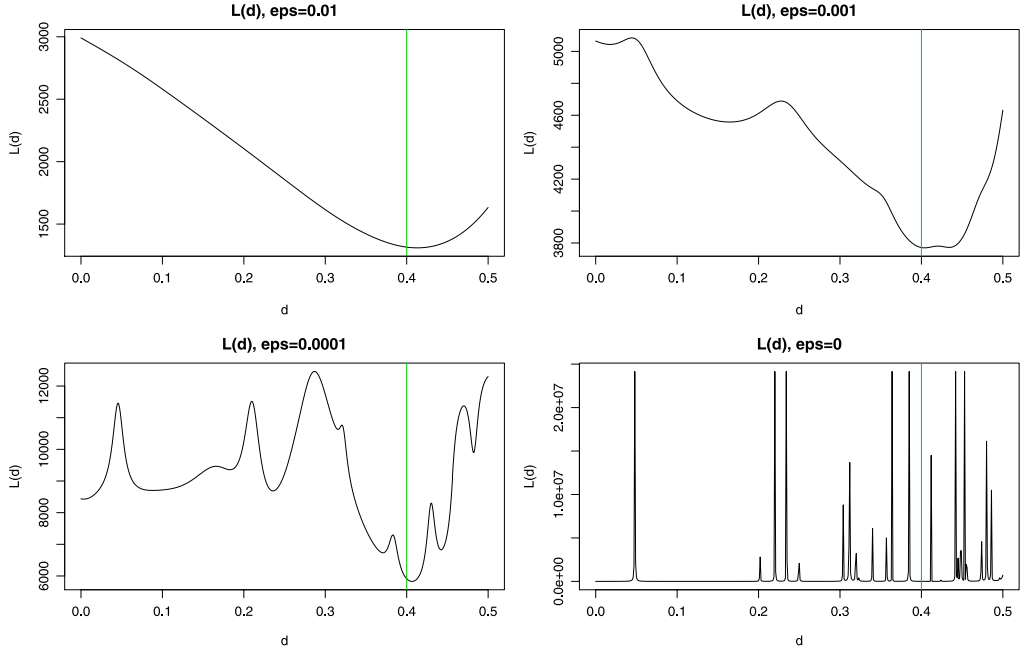

**Figure 1.** For $\epsilon = 0.01, 0.001, 0.0001$ and $0$, the function $\tilde{L}_n$ is plotted as a function of $d$ with fixed $a = 1$ and $c = 0.1$. In each plot, the same path of $X_t$ is used, where the true parameter value is $\theta_0 = (1, 0.4, 0.1)^T$ and $n = 2000$. The vertical line indicates the true value of $d$.

**Remark 7.** If $d > 0$ is close to zero, then the best rate of convergence $n^{\beta/2}$ is close to $n^{1/2}$. However, for strong long memory with $d$ close to $1/2$, the upper bound for $\beta$, given by $1 - 2d$, is very small. Thus, the number of $\sigma_t$'s used for estimation is very small compared to $n$ and the rate of convergence of $\theta_n^{(\beta)}$ is very slow.

**Remark 8.** Though consistency holds for all $\beta \in (0, 1]$, the asymptotic distribution of $\theta_n^{(\beta)}$ for $\beta \geq 1 - 2d$ remains an open problem. The reason for the bound $1 - 2d$ is that, defining

$$d_n := \sqrt{n}(\tilde{L}'_n(\theta) - \bar{\tilde{L}}'_n(\theta_0)),$$

we have

$$E[|d_n|] = \mathrm{O}(n^{\beta/2 + d - 1/2}),$$

which is $o(1)$ for $\beta < 1 - 2d$. For $\beta = 1 - 2d$, the difference is bounded, but it is unclear whether or not it converges to zero.

***Remark 9.*** Alternative estimates of $\theta_0$ could be defined via moment estimation. For instance, empirical estimates of the first three autocovariances of $X_t^2$, $\gamma_{X^2}(0), \gamma_{X^2}(1)$ and $\gamma_{X^2}(2)$, could be used to estimate $\theta_0$ by the method of moments. Limit theorems in Berkes and Horvath (2003) can then be used to show that the resulting estimate is asymptotically normal and the rate of convergence is $n^{1/2-d}$. This is exactly the rate obtained for $\theta_n^{(\beta)}$ at the border $\beta = 1 - 2d$.

## 4. Simulations

We illustrate Theorem 4 by calculating $\theta_n^{(\beta)}$ for simulated LARCH processes with standard normal $\varepsilon_t$ and a parametrization such that (B1) and (B2) hold. The model parameter vector $\theta$ and the constants $\epsilon$ and $\beta$ are chosen as follows:

- Case 1: $d = 0.1$, $a = 1$, $c = 0.2$; $\epsilon = 0.01$, $\beta = 0.799$;
- Case 2: $d = 0.2$, $a = 1$, $c = 0.2$; $\epsilon = 0.01$, $\beta = 0.599$.

To simulate the process $X_t$ via (1) and (2), a pre-sample of length 10 000 is used for initiation. Moreover, the infinite series in (2) is truncated at order 2000. Figures 2a and b show typical sample paths of $X_t$ for the two cases. The corresponding sample autocorrelation functions of $X_t^2$ are given in Figures 2c and d, respectively.

For simplicity, we focus on the estimation of $d$ only. The asymptotic standard deviation given in Theorem 4b (calculated by simulation) is equal to 1.68 in Case 1 and to 1.14 in Case 2. To compare asymptotic with finite-sample results, a small simulation study is carried out as follows. For sample sizes $n = 1000$, 2500, 5000 and 10 000, $N = 1000$ independent samples of the LARCH process are drawn and the estimator $\theta_n^{(\beta)}$ is calculated. Summary statistics of the results are given in Tables 1 (Case 1) and 2 (Case 2). Normal probability plots based on all 1000 simulations are given in Figures 3a–h and 4a–h.

Comparing the results, one can see a strong discrepancy between robust and non-robust estimates of the expected value, standard deviation and skewness of $\theta_n^{(\beta)}$. The robust estimates are close to the asymptotic values obtained from Theorem 4b, already for $n = 1000$. This is not the case for the non-robust estimates. Most extreme are the values of the (non-robust) skewness measure which should converge to zero, but instead seem to be increasing in absolute value. This can be explained as follows. Out of $N = 1000$ simulations, there are a few cases where the algorithm terminated at a solution equal, or very close to, the lower end of the parameter range used in the numerical minimization (see also Remark 5 and Figure 1). As expected from Theorem 4a (and b), the number of cases where this happens decreases with increasing $n$. However, since the variance of estimates in the interior of $\Theta$ tends to zero with increasing $n$, those few estimates that are equal to the *fixed* lower limit of the parameter space become increasingly extreme outliers, compared to the bulk of the simulated data. Indeed, even if $N$ tends to infinity and only one out of $N$ simulations is equal to the lower bound, the empirical skewness will not converge to zero. For this reason, the (non-robust) empirical standard deviation, skewness and normal probability plot are grossly contaminated by the small (and asymptotically negligible) number of simulations where the algorithm did not converge properly. Apart from the robust estimates, we therefore also computed the same empirical non-robust quantities leaving out the ten (out of $N = 1000$) smallest values of $\theta_n^{(\beta)}$. The non-robust estimates are then indeed much closer to the theoretical values, and the normal probability plots indicate convergence (albeit rather slow for $d = 0.2$) to the normal distribution.

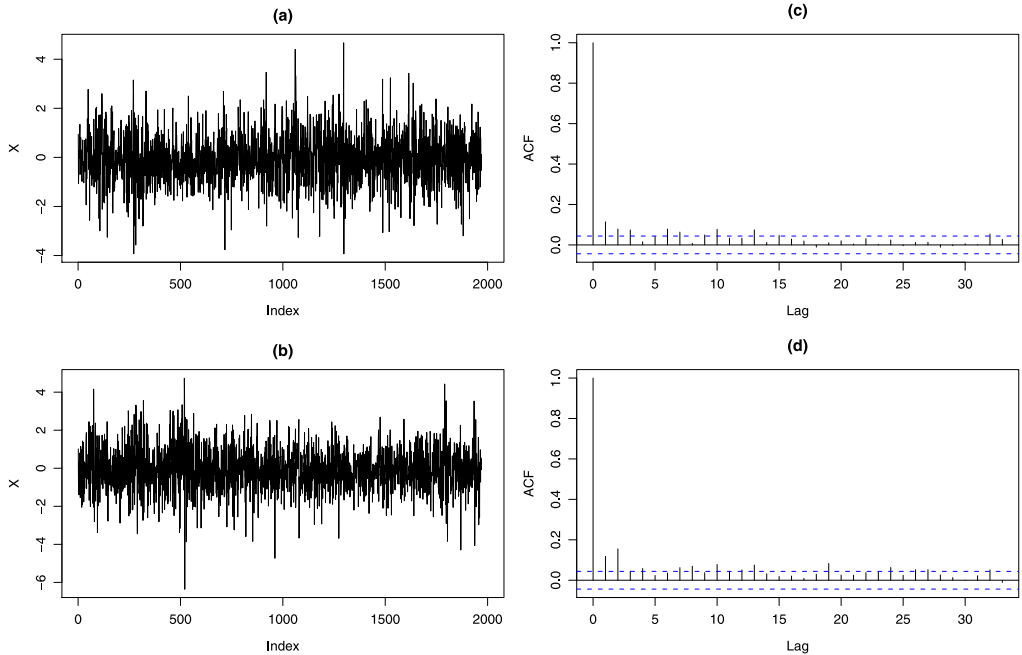

**Figure 2.** Two simulated sample paths of a long-memory LARCH process $X_t$ and the corresponding sample autocorrelation functions of $X_t^2$. The long-memory parameter $d$ is equal to 0.1 in Figures 2a and c, and to 0.2 in Figures 2b and d, respectively.

An additional observation we can make is that convergence to the asymptotic distribution is slower for stronger long memory ($d = 0.2$). The reason is that for $d = 0.2$, the number of terms used in $\tilde{L}_n(\theta)$ is much smaller, namely $\mathrm{O}(n^{0.599})$, as compared to $\mathrm{O}(n^{0.799})$ for $d = 0.1$. More specifically, for $n = 1000$, 2500, 5000 and 10 000, we have $m(n) = 62$, 108, 164 and 248 for $d = 0.2$, whereas for $d = 0.1$, we have $m(n) = 249$, 518, 902 and 1570 for $d = 0.1$.

## 5. Final remarks

We considered parametric estimation for LARCH processes using a modified conditional pseudo-likelihood function. The rate of convergence of the computable version discussed in Section 3.2 depends on the strength of long memory. For short-memory processes ($d \leq 0$), the usual central limit theorem with $\sqrt{n}$-convergence holds. If, on the other hand, $d$ is close to $\frac{1}{2}$, convergence is very slow, so long time series are needed to obtain reliable estimates. In view of the typical range of applications of volatility models, this may not necessarily be a problem. For instance, for high-frequency data in finance, the sample size $n$ is often close to $100\,000$ or more so that the application of $\theta_n^{(\beta)}$ is feasible. How far the best rate $n^{1/2-d}$ may be improved is an open problem. Alternative methods, including Whittle estimation and improved approximations of $\sigma_t$, are the subjects of current research.

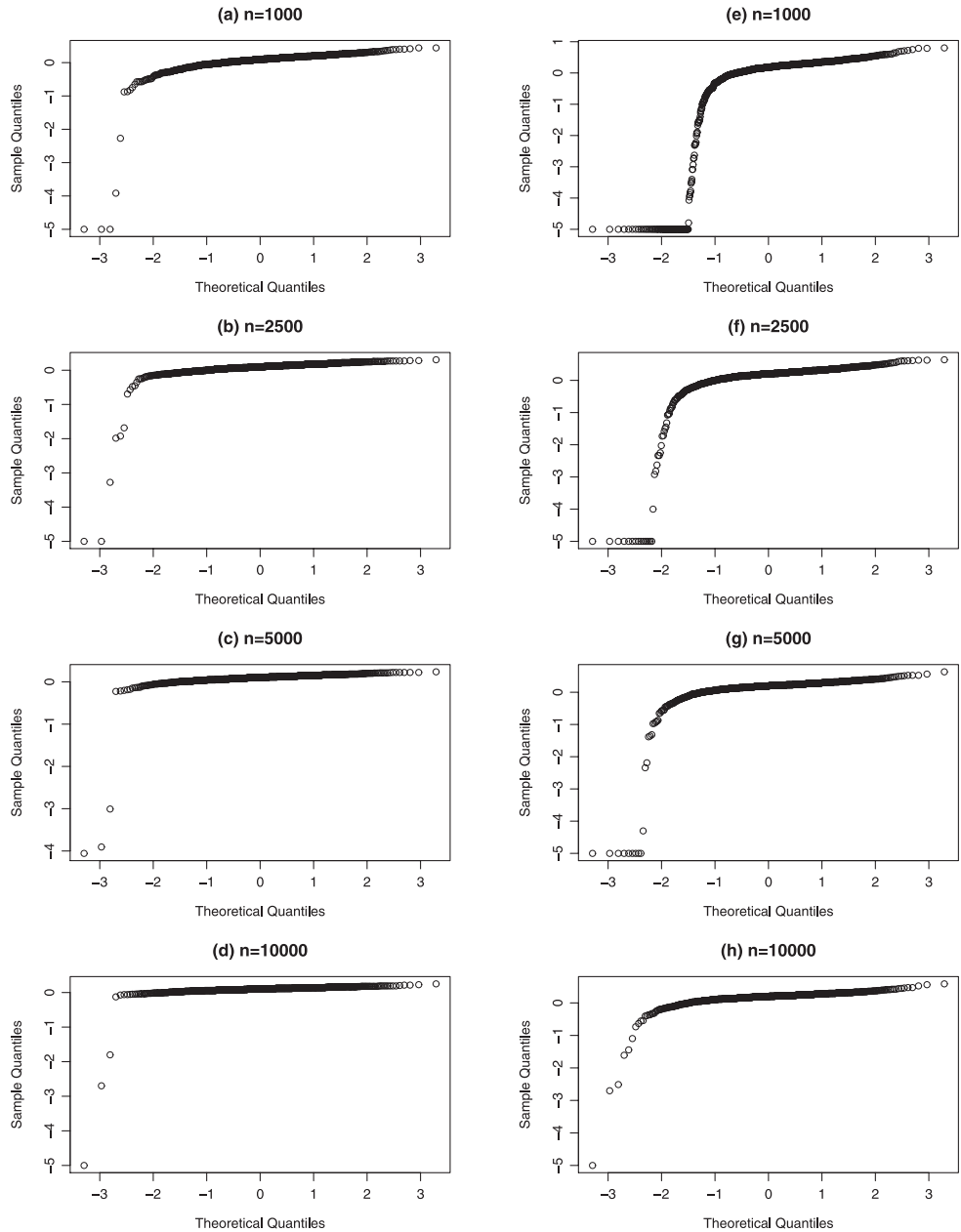

**Figure 3.** Normal probability plots of $N = 1000$ simulated estimates $\theta_n^{(\beta)}$ for Case 1 (Figures 3a–d) and Case 2 (Figures 3e–h).

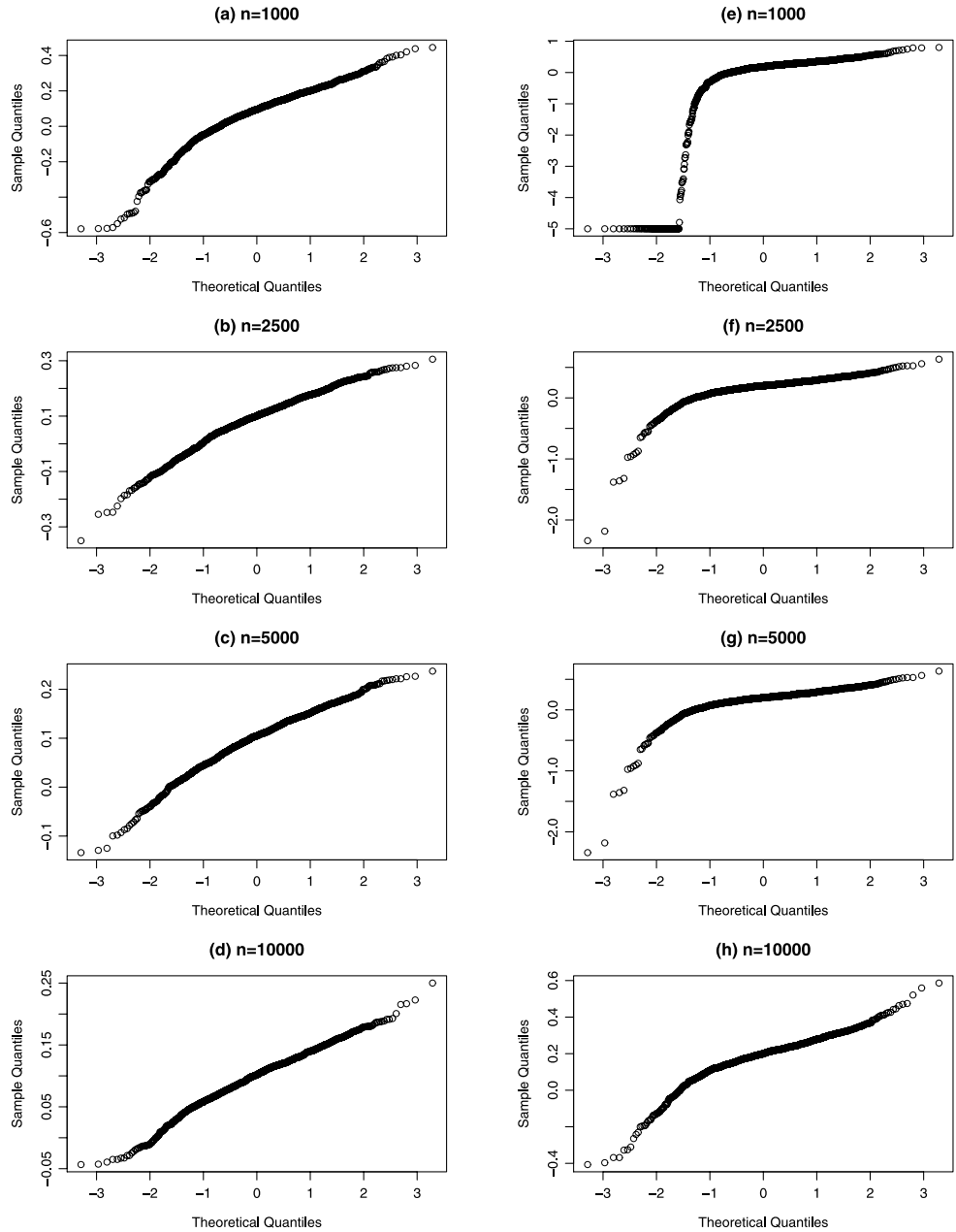

**Figure 4.** Normal probability plots of simulated estimates $\theta_n^{(\beta)}$ for Case 1 (Figures 4a–d) and Case 2 (Figures 4e–h), with ten (out of $N = 1000$) of the lowest points excluded.

**Table 1.** Mean, standard deviation and skewness of $\theta_n^{(\beta)}$ with $\beta = 0.599$, based on $N = 1000$ simulated LARCH processes with long-memory parameter $d = 0.2$ (Case 2). The asymptotic standard deviation from Theorem 4(b) is equal to 1.681. Here, $s$ is the empirical standard deviation, $\tilde{s}$ is the MAD divided by the 75%-percentile of the standard normal distribution and q-skewness is the empirical quartile skewness. In the upper table, all $N = 1000$ simulated values are used; in the lower table, the ten smallest values of $\theta_n^{(\beta)}$ are excluded

| $n$ | 1000 | 2500 | 5000 | 10000 |
|---|---|---|---|---|
| $d = 0.1$: all 1000 simulations | | | | |
| Mean | 0.047 | 0.069 | 0.085 | 0.088 |
| Median | 0.094 | 0.099 | 0.104 | 0.101 |
| $s$ | 0.353 | 0.290 | 0.216 | 0.198 |
| $\tilde{s}$ | 0.121 | 0.082 | 0.054 | 0.041 |
| $n^{\beta/2} s$ | 5.570 | 6.605 | 6.490 | 7.864 |
| $n^{\beta/2} \tilde{s}$ | 1.909 | 1.859 | 1.621 | 1.629 |
| Skewness | $-10.620$ | $-13.161$ | $-16.320$ | $-20.464$ |
| q-skewness | $-0.118$ | $-0.038$ | $-0.093$ | $-0.089$ |
| $d = 0.1$: 10 smallest values of $\hat{d}$ excluded | | | | |
| Mean | 0.072 | 0.091 | 0.098 | 0.098 |
| Median | 0.094 | 0.101 | 0.104 | 0.101 |
| $s$ | 0.150 | 0.090 | 0.057 | 0.043 |
| $\tilde{s}$ | 0.119 | 0.080 | 0.053 | 0.040 |
| $n^{\beta/2} s$ | 2.384 | 2.063 | 1.720 | 1.722 |
| $n^{\beta/2} \tilde{s}$ | 1.882 | 1.822 | 1.595 | 1.604 |
| Skewness | $-1.199$ | $-0.747$ | $-0.684$ | $-0.422$ |
| q-skewness | $-0.103$ | $-0.042$ | $-0.075$ | $-0.079$ |

## Appendix

**Lemma 1.** *Let $(\xi(d, \omega))_{d \in [a,b]}$ be a real-valued separable stochastic process with mean 0 and $E(\xi^2(d)) < \infty$ for all $d \in [a, b]$.*

**Table 2.** Mean, standard deviation and skewness of $\theta_n^{(\beta)}$ with $\beta = 0.599$, based on $N = 1000$ simulated LARCH processes with long-memory parameter $d = 0.2$ (Case 2). The asymptotic standard deviation from Theorem 4(b) is equal to 1.14. Here, $s$ is the empirical standard deviation, $\tilde{s}$ is the MAD divided by the 75%-percentile of the standard normal distribution and q-skewness is the empirical quartile skewness. In the upper table, all $N = 1000$ simulated values are used; in the lower table, the ten smallest values of $\theta_n^{(\beta)}$ are excluded

| $n$ | 1000 | 2500 | 5000 | 10000 |
|---|---|---|---|---|
| | $d = 0.2$: all 1000 simulations | | | |
| Mean | −0.292 | 0.059 | 0.110 | 0.168 |
| Median | 0.181 | 0.201 | 0.198 | 0.199 |
| $s$ | 1.395 | 0.719 | 0.552 | 0.255 |
| $\tilde{s}$ | 0.215 | 0.133 | 0.102 | 0.082 |
| $n^{\beta/2} s$ | 11.041 | 7.489 | 7.079 | 4.030 |
| $n^{\beta/2} \tilde{s}$ | 1.703 | 1.385 | 1.310 | 1.291 |
| Skewness | −2.761 | −5.800 | −7.899 | −11.752 |
| q-skewness | −0.292 | −0.134 | −0.117 | −0.093 |
| | $d = 0.2$: 10 smallest values of $\hat{d}$ excluded | | | |
| Mean | −0.245 | 0.110 | 0.161 | 0.186 |
| Median | 0.184 | 0.202 | 0.199 | 0.200 |
| $s$ | 1.319 | 0.511 | 0.219 | 0.114 |
| $\tilde{s}$ | 0.213 | 0.131 | 0.101 | 0.080 |
| $n^{\beta/2} s$ | 10.437 | 5.319 | 2.810 | 1.800 |
| $n^{\beta/2} \tilde{s}$ | 1.688 | 1.362 | 1.290 | 1.262 |
| Skewness | −2.949 | −6.831 | −4.829 | −1.336 |
| q-skewness | −0.285 | −0.112 | −0.098 | −0.081 |

(a) Denote the covariance function of $\xi$ by $v(d, d') = E(\xi(d)\xi(d'))$. The following then hold:
  (i) If $v(d, d')$ is continuous in $(d, d')$, then $(\xi(d))_{d \in [a,b]}$ is measurable.
  (ii) If $v(d, d')$ is continuously differentiable, then $(\xi(d))_{d \in [a,b]}$ is mean square differentiable, that is, there is a process $(\xi'(d))_{d \in [a,b]}$ with

$$E\left|\frac{1}{h}(\xi(d+h) - \xi(d)) - \xi'(d)\right|^2 \xrightarrow{h \to 0} 0$$

  for all $d \in [a,b]$. Moreover, for almost all $\omega$, $\xi'(\cdot, \omega)$ coincides with the distributional derivative $\partial \xi(\cdot, \omega)/\partial d$.
  (iii) If $v(d, d')$ is $m$ times continuously differentiable, then, for almost all $\omega$, $\xi(\cdot, \omega)$ is $m - 1$ times continuously differentiable.

(b) *If $(\xi(d))_{d\in[a,b]}$ is mean square differentiable with $E(|\xi(d)|^m) < \infty$ and $E(|\xi'(d)|^m) < \infty$ for $m \geq 1$, then*

$$E\left(\sup_{d\in[a,b]} |\xi(d)|^m\right) \leq E(|\xi(a)|^m) + E(|\xi(b)|^m)$$
$$+ m(b-a) \sup_{d\in[a,b]} \{E(|\xi(d)|^m)\}^{(m-1)/m}\{E(|\xi'(d)|^m)\}^{1/m}.$$

**Proof.** (a) is from Kunita (1990), page 40, whereas (iii) is essentially an application of Sobolev's embedding theorem (see, for example, Adams and Fournier (2003)). (b) is an extension of Theorem 3B in Parzen (1965), page 85. □

Throughout this appendix, $K_i$ will denote generic finite constants.

**Lemma 2.** *Let $\theta = (\theta_1, \theta_2, \theta_3)^T$ and suppose that* (A1), (B1), (B2) *and* (S) *hold. Then:*

(a) *under* $(M_3)$, $(M'_p)$ *or* $(M''_p)$, *we have for $k \leq 3$*

$$E\left(\sup_{\theta\in\Theta}\left|\frac{\partial^k \sigma_t(\theta)}{\partial \theta_{i_1}\cdots\partial\theta_{i_k}}\right|^p\right) < \infty,$$

*where $p = 3$ if* $(M_3)$ *holds;*

(b) *under* $(M_3)$, $(M'_p)$ *or* $(M''_p)$, *we have*

$$E\left(\sup_{\theta\in\Theta} |\sigma_t(\theta) - \bar{\sigma}_t(\theta)|^p\right) \to 0 \qquad as\ t \to \infty,$$

*where $p = 3$ if* $(M_3)$ *holds.*

**Proof.** We only give the proof under $(M_3)$. The proof is a combination of Lemma 1b and the combinatorial arguments of Lemmas B.1–B.3 from Giraitis *et al.* (2003). The other cases follow by similar arguments and by using, under $(M'_p)$, Lemma 3.1 of Giraitis *et al.* (2000b) and, under $(M''_p)$, Proposition 2.2 of Giraitis *et al.* (2003), respectively. First, note that

$$\dot{\sigma}_t^{(i_1,\ldots,i_k)}(\theta) = \frac{\partial^k \sigma_t(\theta)}{\partial \theta_{i_1}\cdots\partial\theta_{i_k}}.$$

Moreover, $\sigma_t(\theta) - \bar{\sigma}_t(\theta)$ can be expanded as a Volterra series of the type

$$\Phi_t := \sum_{k=1}^{\infty} \Phi_t^{(k)}$$

with

$$\Phi_t^{(k)} = \sum_{s_k<\cdots<s_1<t} f_{t,1}(t-s_1)f_{t,2}(s_1-s_2)\cdots f_{t,2}(s_{k-1}-s_k)\varepsilon_{s_1}\cdots\varepsilon_{s_k},$$

$f_{t,1}, f_{t,2} \in L^2(\mathbb{Z}_+^0)$, $\|f_{t,1}\|_2 < \infty$ and $\|f_{t,2}\|_2 < 1$. For (a), we set

$$\dot{\sigma}_t^{(i_1,\ldots,i_k)}(\theta) = a\Phi_t + \mathbf{1}_{\{k=1, \theta_1=a\}}$$

with

$$(f_{t,1}(j))_{j\geq 1} = \left(\frac{\partial^k}{\partial \theta_{i_1} \cdots \partial \theta_{i_k}} b_j(\theta)\right)_{j\geq 1}$$

and

$$(f_{t,2}(j))_{j\geq 1} = (b_j(\theta))_{j\geq 1},$$

while for (b),

$$\sigma_t(\theta) - \bar{\sigma}_t(\theta) = a\Phi_1$$

with

$$(f_{t,1}(j))_{j\geq 1} = (b_{j+t}(\theta))_{j\geq 1}$$

and

$$(f_{t,2}(j))_{j\geq 1} = (b_j(\theta))_{j\geq 1}.$$

The proof then follows from the application of Lemma 1b and the following result. A small modification of Lemmas B.1 and B.3 in Giraitis *et al.* (2003) shows that

$$E|\Phi_t|^3 \leq \sum_{k_1, k_2, k_3=1}^{\infty} E[|\Phi_t^{(k_1)} \Phi_t^{(k_2)} \Phi_t^{(k_3)}|]$$

and

$$E[|\Phi_t^{(k_1)} \Phi_t^{(k_2)} \Phi_t^{(k_3)}|] \leq D_{t,1}^3 D_{t,2}^{k_1+k_2+k_3-3},$$

where

$$D_{t,i} = |\mu|_3^{1/3} \|f_{t,i}\|_3 + 3\zeta \|f_{t,i}\|_2$$

and $\zeta$ is defined as in assumption (M$_3$). Hence,

$$E|\Phi_t|^3 \leq \frac{D_{t,1}^3}{(1 - D_{t,2})^3}.$$

Since $\Theta$ is compact, we get in (a) that $D_{t,1} < C_1$ and $D_{t,2} < 1 - C_2$, where the constants $C_1 < \infty$ and $0 < C_2 < 1$ are independent of $\theta$. Furthermore, in (b), $D_{t,1} \to 0$ as $t \to \infty$, uniformly for all $\theta \in \Theta$. Note that $\|f_{t,1}\|_2$ may be greater than 1 and only $\|f_{t,2}\|_2 < 1$ is used. □

**Lemma 3.** *Let assumptions* (A1), (B1), (B2) *and* (S) *hold. Then, under* (M$_3$) *or* (M$'_3$),

$$\sup_{\theta \in \Theta} |L_n(\theta) - L(\theta)| \to 0 \qquad a.s.\ as\ n \to \infty. \tag{3}$$

*If* (M$''_4$) *holds, then*

$$\sup_{\theta \in \Theta} \|L'_n(\theta) - L'(\theta)\| \to 0 \qquad a.s.\ as\ n \to \infty, \tag{4}$$

*where* $L'(\theta) = E(\frac{\partial}{\partial \theta} l_t(\theta))$. *If* (M$'_5$) *holds, then*

$$\sup_{\theta \in \Theta} \|L''_n(\theta) - L''(\theta)\| \to 0 \qquad a.s.\ as\ n \to \infty, \tag{5}$$

*where* $L''(\theta) = E(\frac{\partial^2}{\partial \theta\, \partial \theta'} l_t(\theta))$. *In the three respective cases,* $L(\theta)$ *(resp.* $L'(\theta), L''(\theta)$*) is continuous in* $\theta$.

**Proof.** We first prove (3). From (B1), we have

$$\sup_{\theta \in \Theta} E|l_t(\theta)| \le K(E(X_t^2) + \epsilon) + K \sup_{\theta \in \Theta} E(\sigma_t^2(\theta)) < \infty.$$

Thus $L_n(\theta) \overset{\text{a.s.}}{\to} L(\theta)$ by ergodicity of $X_t^2$ and $\sigma_t(\theta)$ for each individual $\theta \in \Theta$. Uniform convergence follows from a.s. equicontinuity of $(L_n(\theta))_{\theta \in \Theta}$. From the mean value theorem, and the stationarity and ergodicity of $\frac{\partial}{\partial \theta} l_t(\theta)$, it suffices to show that

$$E\left(\sup_{\theta \in \Theta} \left\| \frac{\partial}{\partial \theta} l_t(\theta) \right\| \right) < \infty$$

(see, for example, Andrews (1992)). Since

$$\left\| \frac{\partial}{\partial \theta} l_t(\theta) \right\| \le K_1 |\partial_d \sigma_t(\theta)| X_t^2 + K_2,$$

we get from Hölder's inequality and Lemma 2 that

$$E\left(\sup_{\theta \in \Theta} \left\| \frac{\partial}{\partial \theta} l_t(\theta) \right\| \right) \le K_1 \{E(|X_t|^3)\}^{2/3} \left\{ E\left( \sup_{0 \le d \le d_u} |\partial_d \sigma_t(d)|^3 \right) \right\}^{1/3} + K_2 < \infty.$$

In (4) and (5), pointwise convergence again follows from ergodicity and the particular moment assumption. Uniform convergence is also proved as above. Note that the Hessian matrix of $L_n(\theta)$ is given by

$$L''_n(\theta) = \frac{1}{n} \sum_{t=1}^{n} \frac{\partial^2}{\partial \theta\, \partial \theta'} l_t(\theta),$$

where

$$l_t''(\theta) = \frac{\partial^2}{\partial\theta\,\partial\theta'}l_t(\theta)$$

$$= \frac{4\sigma_t^2(\theta)}{(\sigma_t^2(\theta)+\epsilon)^2}\left(2\frac{X_t^2+\epsilon}{\sigma_t^2(\theta)+\epsilon}-1\right)\frac{\partial}{\partial\theta}\sigma_t(\theta)\left(\frac{\partial}{\partial\theta}\sigma_t(\theta)\right)^T \qquad (6)$$

$$+ \frac{2}{\sigma_t^2(\theta)+\epsilon}\left(1-\frac{X_t^2+\epsilon}{\sigma_t^2(\theta)+\epsilon}\right)\left[\frac{\partial}{\partial\theta}\sigma_t(\theta)\left(\frac{\partial}{\partial\theta}\sigma_t(\theta)\right)^T + \sigma_t(\theta)\frac{\partial^2}{\partial\theta\,\partial\theta'}\sigma_t(\theta)\right].$$

Hence, the matrix norm of $l_t''(\theta)$ is dominated by a linear combination of the terms

$$\sup_{\theta\in\Theta}\left|\frac{\partial}{\partial\theta_i}\sigma_t(\theta)\frac{\partial}{\partial\theta_i}\sigma_t(\theta)X_t^2\right|$$

and

$$\sup_{\theta\in\Theta}\left|\frac{\partial^2}{\partial\theta_i\,\partial\theta_j}\sigma_t(\theta)X_t^2\right|$$

for $i,j \in \{1,2,3\}$. Under (M$_4''$) and Lemma 2, these are bounded in $L^1$ so that (4) follows. Analogously, under (M$_5'$), (5) follows by the $L^1$-boundedness of a similar linear combination also involving the terms

$$\sup_{\theta\in\Theta}\left|\frac{\partial^3}{\partial\theta_i\,\partial\theta_j\,\partial\theta_k}\sigma_t(\theta)X_t^2\right|, \qquad \sup_{\theta\in\Theta}\left|\frac{\partial}{\partial\theta_i}\sigma_t(\theta)\frac{\partial^2}{\partial\theta_j\,\partial\theta_k}\sigma_t(\theta)X_t^2\right|$$

and

$$\sup_{\theta\in\Theta}\left|\frac{\partial}{\partial\theta_i}\sigma_t(\theta)\frac{\partial}{\partial\theta_j}\sigma_t(\theta)\frac{\partial}{\partial\theta_k}\sigma_t(\theta)X_t^2\right|,$$

where $i,j,k \in \{1,2,3\}$, for which Lemma 2 can again be applied. □

**Lemma 4.** *Under* (A1), (B1), (B2) *and* (S), *for every* $\theta \in \Theta\setminus\{\theta_0\}$,

$$L(\theta) > L(\theta_0).$$

**Proof.** From $E(\varepsilon_t^2) = 1$, we get

$$L(\theta) - L(\theta_0) = E\left(\frac{\sigma_t^2+\epsilon}{\sigma_t^2(\theta)+\epsilon} - \ln\left(\frac{\sigma_t^2+\epsilon}{\sigma_t^2(\theta)+\epsilon}\right) - 1\right).$$

Since $x - \ln(x) - 1 > 0$ for $1 \neq x > 0$, we have

$$L(\theta) \geq L(\theta_0)$$

for all $\theta$ and $L(\theta) = L(\theta_0)$ if and only if $\sigma_t^2(\theta) = \sigma_t^2(\theta_0)$ almost surely. Given $\theta$ and $\theta_0$ with $\sigma_t^2(\theta) = \sigma_t^2(\theta_0)$ a.s., we show $\theta = \theta_0$. Thus we define the sets

$$A = \{\omega \in \Omega | \sigma_t(\theta) = \sigma_t(\theta_0)\},$$

$$N_t = \{\omega \in \Omega | \sigma_t \neq 0\}$$

and

$$\bar{A} = A^C \cap \{\omega \in \Omega | \sigma_t^2(\theta) = \sigma_t^2(\theta_0)\}.$$

Note that

$$\bar{A} = \{\omega \in \Omega | \sigma_t(\theta) = -\sigma_t(\theta_0)\}.$$

On $\bar{A} \cap N_{t-1}$, we have

$$a + \sum_{j=1}^{\infty} cj^{d-1} X_{t-j} = -a_0 - \sum_{j=1}^{\infty} c_0 j^{d_0-1} X_{t-j}$$

and hence

$$\varepsilon_{t-1} = -\frac{1}{(c_0+c)\sigma_{t-1}} \left\{ a + a_0 + \sum_{j=2}^{\infty} (cj^{d-1} + c_0 j^{d_0-1}) X_{t-j} \right\}.$$

The right-hand side is measurable w.r.t. $\mathcal{F}_{t-2}$ and hence independent of the left-hand side. Since $\varepsilon_{t-1}$ has a continuous distribution, this is only possible if

$$P(\bar{A} \cap N_{t-1}) = 0.$$

On the sets

$$\bar{A}_k = \bar{A} \bigcap_{i=1}^{k-1} N_{t-i}^C \cap N_{t-k}$$

for $k \geq 2$, repeat the same arguments for $\varepsilon_{t-k}$ to conclude that $P(\bar{A}) = 0$. Note that the set $\{\omega \in \Omega | \exists t_0 : \sigma_t = 0 \text{ for all } t \leq t_0\}$ has probability zero, otherwise equation (2) would not hold. Consequently, with probability one, $\sigma_t(\theta) = \sigma_t(\theta_0)$, that is,

$$a - a_0 = \sum_{j=1}^{\infty} (c_0 j^{d_0-1} - cj^{d-1}) X_{t-j}.$$

Expectation yields $a = a_0$. Finally, considering the variance yields $c_0 j^{d_0-1} = cj^{d-1}$ for all $j \geq 1$. □

**Lemma 5.** *Under* (A1), (B1), (B2), (S) *and* (M$_5$), *the matrices* $G_\epsilon$ *and* $H_\epsilon$ *are positive definite for all* $\theta \in \Theta$.

**Proof.** We only prove that $H_\epsilon$ is positive definite. The proof for $G_\epsilon$ follows by the same arguments. Given $\lambda \in \mathbb{R}^3$, we have to show that

$$\lambda^T H_\epsilon \lambda = E\left(\frac{4\sigma_t^2}{(\sigma_t^2 + \epsilon)^2}\lambda^T \dot\sigma_t \dot\sigma_t^T \lambda\right) = E\left(\frac{4\sigma_t^2}{(\sigma_t^2 + \epsilon)^2}(\lambda^T \dot\sigma_t)^2\right) > 0.$$

Assume that there is a $\lambda = (\lambda_1, \lambda_2, \lambda_3)^T \in \mathbb{R}^3$ such that

$$\frac{4\sigma_t^2}{(\sigma_t^2 + \epsilon)^2}(\lambda^T \dot\sigma)^2 = 0$$

almost surely. Then, on the set $\{\omega \in \Omega | \sigma_t \neq 0\}$, we have

$$\lambda_1 + \sum_{j=2}^{\infty}(\lambda_2 j^{d-1} + \lambda_3 \log(j) j^{d-1}) X_{t-j} = -\lambda_2 \varepsilon_{t-1} \sigma_{t-1}.$$

By arguments similar to those used in the proof of Lemma 4, we then get $\lambda = 0$. $\square$

**Lemma 6.** *Let assumptions* (A1), (B1), (B2) *and* (S) *hold. Then, under* (M$_3$) *or* (M$_3'$),

$$\sup_{\theta \in \Theta} |L_n(\theta) - \bar{L}_n(\theta)| \xrightarrow{L^1} 0 \qquad \text{as } n \to \infty. \tag{7}$$

*If* (M$_4''$) *holds, then*

$$\sup_{\theta \in \Theta} \|L_n'(\theta) - \bar{L}_n'(\theta)\| \xrightarrow{L^1} 0 \qquad \text{as } n \to \infty. \tag{8}$$

*If* (M$_5'$) *holds, then*

$$\sup_{\theta \in \Theta} \|L_n''(\theta) - \bar{L}_n''(\theta)\| \xrightarrow{L^1} 0 \qquad \text{as } n \to \infty. \tag{9}$$

**Proof.** From the mean value theorem applied to $(x^2 + \epsilon)^{-1}$ and $\ln(x + \epsilon)$, and since the derivatives of these functions are bounded, we get

$$\sup_{\theta \in \Theta}|\bar{L}_n(\theta) - L_n(\theta)| \leq \frac{1}{n}\sum_{t=1}^n |X_t^2 + \epsilon| \sup_{\theta \in \Theta}\left|\frac{1}{\bar{\sigma}_t^2(\theta) + \epsilon} - \frac{1}{\sigma_t^2(\theta) + \epsilon}\right|$$

$$+ \frac{1}{n}\sum_{t=1}^n \sup_{\theta \in \Theta}|\ln(\bar{\sigma}_t^2(\theta) + \epsilon) - \ln(\sigma_t^2(\theta) + \epsilon)|$$

$$\leq K\left(\frac{1}{n}\sum_{t=1}^n |X_t^2 + \epsilon| \sup_{\theta \in \Theta}|\bar{\sigma}_t(\theta) - \sigma_t(\theta)| + \frac{1}{n}\sum_{t=1}^n \sup_{\theta \in \Theta}|\bar{\sigma}_t(\theta) - \sigma_t(\theta)|\right).$$

Then, by Lemma 2b, $(M_3)$ or $(M_3')$ implies that

$$E\left(\sup_{\theta\in\Theta}|\bar{\sigma}_t(\theta) - \sigma_t(\theta)|^3\right) \to 0.$$

Together with the Cauchy–Schwarz inequality and Cesaro summability, this proves (7). The other limits, (8) and (9), are proved by means of analogous arguments. □

## Acknowledgements

This work was supported in part by the Deutsche Forschungsgemeinschaft (DFG). We wish to thank the referees for their constructive comments.